\documentclass[11pt]{article}
 \usepackage{epsfig}
 \usepackage{amssymb,amsmath,amsthm,amscd}
\usepackage{latexsym}
\usepackage{showlabels}
\usepackage{graphicx}
\usepackage{color}
\usepackage{subfigure}
 \theoremstyle{theorem}
 \newtheorem{thm}{Theorem}[section]

 \newtheorem{prop}[thm]{Proposition}
 
  \newtheorem{remark}[thm]{Remark}
 \title{On a simple model for describing convection of the rotating fluid: integrability, bifurcations and global dynamics}

 \author{Jia Jiao\footnote{College of Science, Dalian Minzu University, 116600 Dalian, China.}\;
        Shuangling Yang\footnote{School of Mathematics, Sichuan Normal University, Chengdu 610066, China} \;
        Qingjian Zhou\footnote{College of Science, Dalian Minzu University, 116600 Dalian, China.}\;
        Kaiyin Huang\footnote{College of Mathematics, Sichuan University, Chengdu 610065, China.(huangky@scu.edu.cn)}\;
}

\date{\empty}

\begin{document}
 \maketitle

 \begin{abstract}
The Glukhovsky-Dolzhansky (GD) model arises naturally from geophysical science,
which describes rotating fluid convection inside the ellipsoid.
This work aims to provide some new insights into the GD model.
(\emph{i}) We first show that,
under some conditions there are homothetic transformations which
covert the GD model into other similar quadric physical models,
therefore,
our results on the GD model can be naturally applied to the investigation of these models.
(\emph{ii}) We propose a complete classification of Darboux polynomials and exponent factors for the GD model,
which implies that the GD model has no polynomial, rational, or Darboux first integrals. In addition,
some integrable cases of the GD model are also given when the physical parameters are allowed to be non-positive.
(\emph{iii}) The existence of global attractor is proved. The stability and local bifurcations of all co-dimension
one and two are investigated.
Particularly, we show that the GD model undergoes two dynamical transitions as the Rayleigh number increases.
(\emph{iv}) To understand the asymptotic behavior of the orbits for the GD model,
we use the Poincar\'{e} compactification method to study its dynamical behavior at infinity. More precisely, we prove that the phase portraits of the GD model at infinity consist of an infinite sequence of periodic solutions and two heteroclinic loops. Our results  may help us better understand the complex
and rich dynamics of rotating fluid convection.
 \end{abstract}

 \newpage
 
 \section{Introduction}

The motion of rotating fluids has gained interesting attention in the past several decades since the rotation strongly
affects the dynamical features of fluids, compared with non-rotating systems \cite{gr1,gr2}. In addition, the study of dynamics of rotating fluids is of crucial importance for
understanding many phenomena in atmosphere, oceanography, planetary physics and astrophysics, see \cite{gr3,gr4}  and references therein.
However, due to the nonlinear nature, our present knowledge on rotating fluids is still poor and
 the previous results are mainly limited to
 the study of fluid motions in rotating, closed containers filled
with liquids in some simple geometries such as annular channels, circular cylinder,  sphere or spherical shell or ellipsoid.

In 1980,
Glukhovsky and Dolzhansky \cite{model} studied the convection of viscous incompressible fluid motion inside the ellipsoid
\begin{equation*}\label{I0}
(\frac{x_1}{a_1})^{2}+(\frac{x_2}{a_2})^{2}+(\frac{x_3}{a_3})^{2}=1,~~~~a_1>a_2>a_3>0,
\end{equation*}
under the condition of stationary inhomogeneous external heating.
They assumed that the ellipsoid rotates with the constant velocity $\Omega_0$ along the axis $a_3$,
and the axis $a_3$ has a constant angle $\alpha$ with the gravity vector $g$.
They also assumed that the temperature difference is generated along the axis $a_1$ and its gradient has a constant value $q_0$.
Denote by $\lambda,\mu,\beta$ the coefficients of viscosity,
heat conduction,
and volume expansion,
respectively.
Then they proposed the following system of ordinary differential equations
\begin{equation}\label{I1}
\begin{cases}
\dot{x}=Ayz+Cz-\sigma x:=P(x,y,z)\\
\dot{y}=-xz+R_a-y:=Q(x,y,z)\\
\dot{z}=-z+xy:=R(x,y,z)
\end{cases}
\end{equation}
as a three-mode model to describe the convection of viscous incompressible fluid motion inside the ellipsoid.
Here
\begin{align*}
&\sigma=\frac{\lambda}{\mu},~~~T_{\alpha}=\frac{\Omega_{0}^{2}}{\lambda^{2}},
~~~R_{a}=\frac{g\beta a_{3}q_{0}}{2a_{1}a_{2}\lambda\mu},\\
&A=\frac{a_{1}^{2}-a_{2}^{2}}{a_{1}^{2}+a_{2}^{2}}\cos^{2}\alpha T_{\alpha}^{-1},~~
C=\frac{a_{1}^{2}-a_{2}^{2}}{a_{1}^{2}+a_{2}^{2}}\cos^{2}\alpha T_{\alpha}^{-1},\\
&x(t)=\mu^{-1}\omega_{3}(t),~~y(t)=\frac{g\beta a_{3}}{2a_{1}a_{2}\lambda \mu}q_{1}(t),~~z(t)=\frac{g\beta a_{3}}{2a_{1}a_{2}\lambda \mu}q_{2}(t),
\end{align*}
$\omega_3(t)$ is the projection of the vector of angular velocity on the axis $a_3$,
$q_{1}(t)$ and $q_{2}(t)$ are the projections of temperature gradients on the axes $a_1$ and $a_2$.
The parameters $\sigma$, $T_{\alpha}$ and $R_{a}$ are the Prandtl,
Taylor and Rayleigh numbers, respectively.
The parameters $A,C,\sigma$ and $Ra$ of system (\ref{I1}) are positive due to their physical meanings.

The nonlinear differential system (\ref{I1}) is now called the Glukhovsky-Dolzhansky (GD) model in the references.
It is significantly different from the classical Lorenz system since the Lorenz system is associated with the flow of the two-dimensional convection but the GD model is associated with the flow of the three-dimensional convection.
In addition,
using the linear transformation \cite{model,gb0}
$$
x\rightarrow x,~y\rightarrow R_aC-\frac{\sigma}{AR_a+C}z,~z\rightarrow \frac{\sigma}{AR_a+C}y,
$$
the GD model is converted into a Lorenz-like form
\begin{equation}\label{lorenzlike}
\begin{cases}
\dot{x}=\sigma(y-x)-ayz\\
\dot{y}=rx-y-xz\\
\dot{z}=-z+xy
\end{cases}
\end{equation}
with
$$
a=\frac{A\sigma^2}{(AR_a+C)^2},~~r=\frac{R_a(AR_a+C)}{\sigma}.
$$
Compared with the classical Lorenz model,
the GD model has an additional nonlinear term,
which has a significant impact on the integrability structures and dynamical features of the GD model.

In recent years,
the GD model has been intensively investigated,
particularly by numerical simulation,
see \cite{gb0,gb1,gb2,gb3} and the references therein.
For instance,
the GD model has a $3$-dimensional strange chaotic attractor
as shown in Figures 1 and 2.
However,
this system has never been studied about its integrability,
local bifurcations or global dynamics.
The main purpose of this work is to cover these gaps and to provide some new insights into the GD model.
Our results show that,
in spite of its simple form,
the GD model admits rich dynamics.

Our main results are summarized as follows:
\begin{itemize}
\item We provide some linear rescaling of time and coordinate to show the GD model can be transformed into the
Rabinovich system or 3D forced-damped system or \emph{D2} vector field only when the parameters satisfy some conditions:
$C=-2AR_a$ or $\sigma=1$ or $R_a=C=0$ respectively.

\item A \emph{complete} classification of invariant algebraic surfaces
and exponent factors is given,
which yields that the GD model is not Darboux integrable for any value of the positive real parameters $A,\sigma,C,R_a$.
However, we also provide some integrable cases of this model when parameters are allowed to be non-positive.

\item Due to the dissipative structure of the GD model,
the existence of its global compact attractive set is proved.
All local bifurcations of co-dimension one and two at equilibrium are investigated.
We show that the GD model undergoes a cusp bifurcation as the Rayleigh number $R_a$
crosses a threshold $R_a^*$.
Moreover, if a nonlinear inequality for $A,\sigma,C$ holds,
the GD model undergoes a Hopf bifurcation as $R_a$ crosses a second threshold $R_a^{**}$.

\item Using the Poincar\'{e} compactification technique,
we give a complete description of the dynamical behavior of the GD model on the sphere at infinity.

\end{itemize}

The paper is organized as follows.
In section 2,
we discuss the relationship between the GD model and some quadric systems.
In section 3,
we study the integrability of the GD model in the framework of Darboux integrability theory.
In section 4,
we focus on the stability and bifurcations of the GD model,
especially the effect of Rayleigh number on the GD model.
In section 5,
the dynamical behaviors of the GD model at infinity are explored by using the classical Poincar\'{e}
compactification for the three-dimensional polynomial vector fields.
Some discussions on the physical meaning of our results are present in the last section.

\section{Homothetic transformation between the GD model and other quadric systems}

There exist some quadric systems e.g. the Rabinovich system \cite{si1,si2},
3D forced-damped system \cite{si3} and \emph{D2} vector field,
which are similar to the GD model.
These systems admit a common structure:
the original is the equilibrium point,
symmetry with respect to a certain coordinate axis,
only three cross nonlinear terms $\{yz,xz,xy\}$.

A key and natural problem is to identify the relationship between the GD model and these systems.
If there exists a homothetic scaling in time and state variables to convert one to the other,
then we can easily get any dynamical object (e.g. equilibrium, periodic orbit,
homoclinic orbit, heteroclinic orbit, Silnikov chaos, attractor set, or chaos synchronization, etc.) from one to the other.

In what follows,
we provide a relationship between the Glukhovsky-Dolzhansky system (\ref{I1}) and other similar quadric
systems including the Rabinovich system,
the 3D forced-damped system and the \emph{D2} vector field \cite{si4}.

The Rabinovich system \cite{si1,si2} is a three-wave interaction model and is given by
\begin{equation}\label{rr}
\begin{cases}
\dot{X}=h Y-v_{1} X+YZ \\
\dot{Y}=h X-v_{2} Y-XZ \\
\dot{Z}=-v_{3} Z+XY,
\end{cases}
\end{equation}
where $h,~v_1,~v_2,~v_3$ are parameters.
When $C=-2ARa$,
by a family of linear scaling $\mathcal{S}_{\alpha}: (x,y,z,t)\rightarrow (X,Y,Z,T)$
$$
x=\frac{X}{\alpha},~~y=Ra-\frac{Z}{\alpha\sqrt{-A}},~~z=\frac{Y}{\alpha\sqrt{-A}},~~t=\alpha T,
$$
we transform the Glukhovsky-Dolzhansky system (\ref{I1}) into
\begin{equation}\label{rr2}
\begin{cases}
\frac{dX}{dT}=\alpha\sqrt{-A}Ra Y-\alpha\sigma X+YZ \\
\frac{dY}{dT}=\alpha\sqrt{-A}Ra X-\alpha Y-XZ \\
\frac{dZ}{dT}=-\alpha Z+XY.
\end{cases}
\end{equation}
Clearly,
system (\ref{rr2}) corresponds to the Rabinovich system (\ref{rr}) with parameters $(h,v_1,v_2,v_3)=(\alpha\sqrt{-A}Ra, \alpha\sigma,\alpha,\alpha)$.
Hence, if $C=-2ARa$,
the Glukhovsky-Dolzhansky system (\ref{I1})
is equivalent to the Rabinovich system (\ref{rr}) in the particular case of the parameter region $\{(h,v_1.v_2.v_3)|v_2=v_3\}$,
that is,
two systems are homothetic copies.
In addition,
the linear scaling $\mathcal{S}_{\alpha}$ forces parameter $A$ to be negative,
which cannot hold in the physical region of parameter $A$.

The 3D forced-damped system \cite{si3} arises in mechanical,
electrical and fluid dynamical contexts and is given by
\begin{equation}\label{3ff}
\begin{cases}
\dot{X}=-aX+Y+YZ \\
\dot{Y}=X-aY+bXZ \\
\dot{Z}=cZ-bXY,
\end{cases}
\end{equation}
where $a,b,c$ are parameters.
Compared with the Glukhovsky-Dolzhansky system (\ref{I1}),
the 3D forced-damped system (\ref{3ff})
has less parameters.
When $\sigma=1$, the linear scaling
\begin{equation*}
x=\frac{(AR_a+C)^{3/2}}{A\sqrt{R_a}}X,~~y=Ra+\frac{(AR_a+C)Z}{A},
\end{equation*}
\begin{equation*}
z=\frac{(AR_a+C)Y}{A},~~t=\frac{T}{\sqrt{(AR_a+C)R_a}}
\end{equation*}
transforms the Glukhovsky-Dolzhansky system (\ref{I1}) into
\begin{equation}\label{3ff2}
\begin{cases}
\frac{dX}{dT}=-\frac{1}{\sqrt{(AR_a+C)R_a}}X+Y+YZ \\
\frac{dY}{dT}=X-\frac{1}{\sqrt{(AR_a+C)R_a}}Y+\frac{AR_a+C}{{AR_a}}XZ \\
\frac{dZ}{dT}=-\frac{1}{\sqrt{(AR_a+C)R_a}}Z-\frac{AR_a+C}{{AR_a}}XY,
\end{cases}
\end{equation}
which corresponds to the 3D forced-damped system  (\ref{3ff}) with parameters
$$
(a,b,c)=(\frac{1}{\sqrt{(AR_a+C)R_a}}, \frac{AR_a+C}{{AR_a}},-\frac{1}{\sqrt{(AR_a+C)R_a}}).
$$
Therefore, when $\sigma=1$, the Glukhovsky-Dolzhansky system (\ref{I1})
is equivalent to the 3D forced-damped system (\ref{3ff}) in the particular case of the parameter region $\{(a,b,c)|a+c=0\}$.

Consider the vector field equivalent under the \emph{D2} symmetry group,
called the \emph{D2} vector field \cite{si4}
\begin{equation}\label{DD20}
\begin{cases}
\dot{X}=aX+YZ \\
\dot{Y}=bY+XZ \\
\dot{Z}=Z-XY,
\end{cases}
\end{equation}
where $a,b$ are parameters.
Similar to above, when $R_a=C=0$, we make a linear scaling
$$
x=X,~~y=\frac{Z}{\sqrt{A}},~~z=-\frac{Y}{\sqrt{A}},~~t=-T
$$
and transform (\ref{I1}) into
\begin{equation}\label{DD22}
\begin{cases}
\frac{dX}{dT}=\sigma X+YZ \\
\frac{dY}{dT}=Y+XZ \\
\frac{dZ}{dT}=Z-XY,
\end{cases}
\end{equation}
which corresponds to the \emph{D2} vector field (\ref{DD20}) with parameters $(a,b)=(\sigma,1)$.
Then, when $R_a=C=0$,
the Glukhovsky-Dolzhansky system (\ref{I1})
is equivalent to the \emph{D2} vector field (\ref{DD20}) in the particular case of the parameter region $\{(a,b)|b=1\}$.

Finally,
we mention that based on the above discussion one can get more information about the Rabinovich system,
3D forced-damped system and \emph{D2} vector field from the results obtained in the papers devoted to the study of the GD model.

\section{Integrability analysis of the GD model}
Generally,
a system of differential equations is integrable if it possesses a sufficient number of first integrals
(and/or other tensor invariants) such that we can solve this system explicitly.
Hence we could obtain its global information and understand its topological structure \cite{dd2,dd3}.
Furthermore,
non-integrability of the system also seems necessary for better understanding of the complex phenomenon \cite{in1,in3}.
However,
to study the integrability of a given system is not an easy task,
since there is no any effective approach to determine the existence or nonexistence of first integrals.

In this paper,
we aim to study the integrability of the GD model in the framework of Darboux integrability theory.
Darboux integrability theory plays an important role in the integrability of
the polynomial differential systems \cite{dar1,dar2,dar3,dar4,dar5},
which helps us find first integrals by knowing a sufficient number of algebraic
invariant surfaces (the Darboux polynomials) and of exponential factors,
see \cite{d1,d2,d3,d4,d5} for instance.
Moreover,
it can also help us make a more precise analysis of the global dynamics of the considered system topologically \cite{dd2,dd1}.

We first recall some basic definitions.
Let $\mathbb{R}[x,y,z]$ be the ring of the real polynomials in the variables $x$, $y$ and $z$.
We say that $f(x,y,z)\in \mathbb{R}[x,y,z]$ is a \emph{Darboux polynomial} of system (\ref{I1}) if it satisfies
\begin{equation}\label{P1}
\frac{\partial f}{\partial x}P+\frac{\partial f}{\partial y}Q+\frac{\partial f}{\partial z}R=Kf,
\end{equation}
for some polynomial $K$, called the cofactor of $f(x,y,z)$.
If $f(x,y,z)$ is a Darboux polynomial,
then the surface $f(x,y,z)=0$ is an invariant manifold of system (\ref{I1}).
Particularly,
if $K=0$, $f(x,y,z)$ satisfies the following equation
\begin{equation}\label{P2}
\frac{\partial f}{\partial x}P+\frac{\partial f}{\partial y}Q+\frac{\partial f}{\partial z}R=0,
\end{equation}
then polynomial $f(x,y,z)$ is called a \emph{polynomial first integral} of system (\ref{I1}) .
Let $g,h\in \mathbb{R}[x,y,z]$ be coprime.
We say that a nonconstant function $E=\exp(g/h)$ is an \emph{exponential factor} of system (\ref{I1})
if $E$ satisfies
\begin{equation*}\label{P3}
\frac{\partial E}{\partial x}P+\frac{\partial E}{\partial y}Q+\frac{\partial E}{\partial z}R=LE,
\end{equation*}
for some polynomial $L\in\mathbb{R}[x,y,z]$ with the degree at most one, called the cofactor of $E$.
A first integral $G$ of system (\ref{I1}) is called \emph{Darboux type} if it is a first integral of the form
\begin{equation*}\label{P4}
G=f_{1}^{\lambda_{1}}\cdots f_{p}^{\lambda_{p}}E_{1}^{\mu_{1}}\cdots E_{q}^{\mu_{q}},
\end{equation*}
where $f_{1},\cdots,f_{p}$ are Darboux polynomials,
$E_{1},\cdots,E_{q}$ are exponential factors and $\lambda_{i},\mu_{j}\in R$,
for $i=1,\cdots,p$ and $j=1,\cdots,q$.

\begin{prop}\label{T1}
System (\ref{I1})  has no polynomial first integrals.
\end{prop}
\begin{proof}Suppose
\begin{align}\label{poly}
f(x,y,z)=\sum_{i=0}^{n}f_{i}(x,y,z)
\end{align}
is a polynomial first integral of system (\ref{I1}), where $f_i$ are the homogeneous polynomials of degree $i$ and $f_{n}\neq0$.
Firstly, substituting $(\ref{poly})$ into (\ref{P2}) and identifying the homogeneous components of degree $n+1$, we get
\begin{equation}\label{po}
Ayz\frac{\partial f_{n}}{\partial x}-xz\frac{\partial f_n}{\partial y}+xy\frac{\partial f_n}{\partial z}=0.
\end{equation}
The characteristic equations associated with (\ref{po}) are
\begin{equation*}\label{D5}
\frac{dx}{dy}=\frac{Ayz}{-xz},~~\frac{dz}{dy}=\frac{xy}{-xz}.
\end{equation*}
Their general solutions are
\begin{equation*}\label{D6}
x^{2}+Ay^{2}=c_{1},~~y^{2}+z^{2}=c_{2},
\end{equation*}
where $c_{1}$ and $c_{2}$ are arbitrary constants.
We make the change of variables
\begin{equation*}\label{D7}
u=x^{2}+Ay^{2},~~w=y,~~v=y^{2}+z^{2}.
\end{equation*}
Correspondingly, the inverse transformation is
\begin{equation}\label{D8}
x=\pm\sqrt{u-Aw^{2}},~~y=w,~~z=\pm\sqrt{v-w^{2}}.
\end{equation}
Without loss of generality, we set
\begin{equation}\label{D9}
x=\sqrt{u-Aw^{2}},~~y=w,~~z=-\sqrt{v-w^{2}},
\end{equation}
and transform (\ref{po}) into
\begin{align}\label{D12}
\frac{d\bar{f}_n}{d w}=0,
\end{align}
where $\bar{f}_n(u,v,w)=f_n(x,y,z)$. In the following, unless otherwise specified,
we always denote the function $R(x,y,z)$ by $\bar{R}(u,v,w)$.
Hence we obtain
\begin{equation*}\label{D14}
f_n(x,y,z)=f_{2m}(x,y,z)=\sum_{i=0}^{m}a_i^m(x^{2}+Ay^{2})^{m-i}(y^{2}+z^{2})^{i},~~a_i^m\in\mathbb{R},
\end{equation*}
where $n=2m$ must be an even number.
Secondly, substituting $(\ref{poly})$ into (\ref{P2}) and identifying the homogeneous components of degree $n$ yields
\begin{align}
Ayz\frac{\partial f_{n-1}}{\partial x}&-xz\frac{\partial f_{n-1}}{\partial y}+xy\frac{\partial f_{n-1}}{\partial z}
=(\sigma x-Cz)\frac{\partial f_{n}}{\partial x}+y\frac{\partial f_{n}}{\partial y}+z\frac{\partial f_{n}}{\partial z}\label{assd}\\
=&\sum_{i=0}^{m}a_{i}^{m}2\sigma(m-i)(x^2+Ay^2)^{m-i-1}(y^2+z^2)^ix^2\notag\\
&-\sum_{i=0}^{m}a_{i}^{m}2C(m-i)(x^2+Ay^2)^{m-i-1}(y^2+z^2)^ixz\notag\\
&+\sum_{i=0}^{m}a_{i}^{m}2A(m-i)(x^2+Ay^2)^{m-i-1}(y^2+z^2)^{i}y^{2}\notag\\
&+\sum_{i=0}^{m}a_{i}^{m}2i(x^2+Ay^2)^{m-i}(y^2+z^2)^{i-1}y^{2}\notag\\
&+\sum_{i=0}^{m}a_{i}^{m}2i(x^2+Ay^2)^{m-i}(y^2+z^2)^{i-1}z^{2}.\notag
\end{align}
Using the transformations (\ref{D9}) again, the above equation becomes
\begin{align*}
\sqrt{u-Aw^2}\sqrt{v-w^2}\frac{d\bar{f}_{2m-1}}{dw}=&\sum_{i=0}^{m}a_{i}^{m}2\sigma(m-i)u^{m-i-1}v^i(u-Aw^2)\\
&+\sum_{i=0}^{m}a_{i}^{m}2C(m-i)u^{m-i-1}v^i\sqrt{u-Aw^2}\sqrt{v-w^2}\\
&+\sum_{i=0}^{m}a_{i}^{m}2A(m-i)u^{m-i-1}v^{i}w^{2}+\sum_{i=0}^{m}a_{i}^{m}2iu^{m-i}v^{i-1}w^{2}\\
&+\sum_{i=0}^{m}a_{i}^{m}2iu^{m-i}v^{i-1}(v-w^2).
\end{align*}
We obtain
\begin{align*}
\bar{f}_{2m-1}
=&\sum_{i=0}^{m}a_{i}^{m}2\sigma(m-i)u^{m-i-1}v^i\int\frac{\sqrt{u-Aw^{2}}}{\sqrt{v-w^2}}dw\\
&+\sum_{i=0}^{m-1}2[a_{i}^{m}A(m-i)+a_{i+1}^{m}(i+1)]u^{m-i-1}v^{i}\int\frac{w^{2}dw}{\sqrt{u-Aw^{2}}\sqrt{v-w^2}}\\
&+\sum_{i=0}^{m}a_{i}^{m}2iv^{i-1}\int\frac{\sqrt{v-w^{2}}}{\sqrt{u-Aw^2}}dw\\
&+\sum_{i=0}^{m}a_{i}^{m}2C(m-i)u^{m-i-1}v^iw+A_{n-1}(u,v),
\end{align*}
where $A_{n-1}(u,v)$  is an arbitrary smooth function in $u$ and $v$.
It is easy to check
\begin{eqnarray}\label{D15}
\int\frac{w^2dw}{\sqrt{u-Aw^2}\sqrt{v-w^2}}
=-\int\frac{\sqrt{u-Aw^2}}{\sqrt{v-w^2}}d w
+u\int\frac{d w}{\sqrt{u-Aw^2}\sqrt{v-w^2}}.
\end{eqnarray}
Since
\begin{align*}
\int\frac{\sqrt{u-Aw^2}}{\sqrt{v-w^2}}dw,~~~~~
\int\frac{ d w}{\sqrt{u-Aw^2}\sqrt{v-w^2}}d w
\end{align*}
are elliptic integrals of the second and first kinds, respectively,
in order that $f_{2m-1}(x,y,z) =\bar f_{2m-1}(u,v,w)$ is a homogeneous polynomial of degree $2m-1$, we must have
\begin{equation}\label{D16}
\begin{cases}
a_{i}^{m}2\sigma(m-i)=0,~~~i=0,1,\cdots,m\\
a_{i}^{m}2C(m-i)=0, ~~~i=0,1,\cdots,m\\
a_{i}^{m}2i=0, ~~~i=0,1,\cdots,m\\
a_{i}^{m}A(m-i)+a_{i+1}^{m}(i+1)=0, ~~~i=0,1,\cdots,m-1.
\end{cases}
\end{equation}
By (\ref{D16}) and $A,C,\sigma$ being positive,
we obtain $a_{i}^{m}=0, i=0,\cdots,m$.
This leads to a contradiction.
So system (\ref{I1})  has no polynomial first integrals.
\end{proof}

\begin{prop}\label{T2}
System (\ref{I1}) has no Darboux polynomial with  nonzero cofactors.
\end{prop}
\begin{proof}
Suppose
\begin{align}\label{aa1}
f(x,y,z)=\sum_{i=0}^{n}f_{i}(x,y,z)
\end{align}
is a Darboux polynomial of the system (\ref{I1}) with a non-cofactor $K(x,y,z)$, where $f_{i}$ is a homogeneous polynomial of degree $i$ for $i=0,1,\cdots,n$. Comparing the degree on both sides of (\ref{P1}) yields $\deg K\leq 1$.
Without loss of generality,
we can assume that the cofactor is of the form
\begin{equation}\label{DD2}
K(x,y,z)=k_1x+k_2y+k_3z+k_0,~~k_i\in\mathbb{R},~~i=0,1,2,3.
\end{equation}
Substituting (\ref{aa1}) and $(\ref{DD2})$ into  (\ref{P1}) and identifying the terms of the same degree,
we obtain
\begin{align}\label{DD3}
Ayz\frac{\partial f_{n}}{\partial x}-xz\frac{\partial f_n}{\partial y}+xy\frac{\partial f_n}{\partial z}&=(k_1x+k_2y+k_3z)f_n,\\\nonumber
Ayz\frac{\partial f_{n-1}}{\partial x}-xz\frac{\partial f_{n-1}}{\partial y}+xy\frac{\partial f_{n-1}}{\partial z}=&(k_1x+k_2y+k_3z)f_{n-1}+(\sigma x-Cz)\frac{\partial f_{n}}{\partial x}\\\nonumber
&+y\frac{\partial f_{n}}{\partial y}
+z\frac{\partial f_{n}}{\partial z}+k_0f_{n},\\\label{DD4}
Ayz\frac{\partial f_{i}}{\partial x}-xz\frac{\partial f_i}{\partial y}+xy\frac{\partial f_i}{\partial z}=&(k_1x+k_2y+k_3z)f_i+(\sigma x-Cz)\frac{\partial f_{i+1}}{\partial x}\\\nonumber
&+y\frac{\partial f_{i+1}}{\partial y}
+z\frac{\partial f_{i+1}}{\partial z}-Ra\frac{\partial f_{i+2}}{\partial y}+k_0f_{i+1}, \\
&~~~i=n-2,\ldots,0.\nonumber
\end{align}

We claim that the cofactor is a constant, i.e., $k_1=k_2=k_3=0$. Indeed, under the change of (\ref{D8}),
we can transform (\ref{DD3}) into an ordinary differential equation if we fixe $u$ and $v$
\begin{equation*}\label{DD6}
-(\pm\sqrt{u-Aw^{2}})(\pm\sqrt{v-w^{2}})\frac{d\bar{f_n}}{dw}
=[k_1(\pm\sqrt{u-Aw^{2}})+k_2w+k_3(\pm\sqrt{v-w^{2}})]\bar{f}_n,
\end{equation*}
where $\bar{f}_n(u,v,w)=f_n(x,y,z)$.
In the following proof, we only consider the case of $xz<0$.
For the case of $xz>0$,
the proof is similar.
Solving the last equation we find that for $xz<0$,
\begin{align*}
\bar{f}_n=&\bar{A}(u,v)\left|\frac{2\sqrt{A}\sqrt{(u-Aw^2)(v-w^2)}+2Aw^2-(u+Av)}{2\sqrt{A}}
\right|^{-\frac{k_2}{2\sqrt{A}}}\\
&\exp\left(-k_1\arcsin\frac{w}{\sqrt{v}}\right)\exp\left(-\frac{k_3}{\sqrt{A}}\arcsin\frac{w}{\sqrt{u}}\right).
\end{align*}
In order for $f_n(x,y,z)=\bar{f_n}(u,v,w)$ to be a homogeneous polynomial of degree $n$ in $x,y,z$,
we have $k_1=k_3=0$ and the function $\bar{A}$ is a homogeneous polynomial in $x^2+Ay^2$ and $y^2+z^2$.
Then
\begin{eqnarray*}
f_n=\bar{A}(x^2+Ay^2,y^2+z^2)(x+\sqrt{A}z)^{\frac{k_2}{\sqrt{A}}}.
\end{eqnarray*}
Therefor $f$ is a Darboux polynomial of degree $n(=2m+\frac{k_2}{\sqrt{A}})$
with the cofactor $K=k_2y+k_0$.
Set $\bar{k}_2={k_2}/{\sqrt{A}}$, we get
$
f=\sum_{i=0}^{2m+\bar{k}_2}f_i
$
and
\begin{eqnarray*}
f_{2m+\bar{k}_2}=(x+\sqrt{A}z)^{\bar{k}_2}\sum_{i=0}^{m}a_{i}^{m}(x^2+Ay^2)^{m-i}(y^2+z^2)^i.
\end{eqnarray*}
Substituting $f_{2m+\bar{k}_2}$ into (\ref{DD4}) and performing some calculations,
we obtain
\begin{align*}
Ayz&\frac{\partial f_{2m+\bar{k}_2-1}}{\partial x}-xz\frac{\partial f_{2m+\bar{k}_2-1}}{\partial y}+xy\frac{\partial f_{2m+\bar{k}_2-1}}{\partial z}-k_2yf_{2m+\bar{k}_2-1}\\
&=(x+\sqrt{A}z)^{\bar{k}_2}\sum_{i=0}^{m}a_{i}^{m}[2\sigma(m-i)+k_{0}+2i](x^2+Ay^2)^{m-i}(y^2+z^2)^i\\
&+(x+\sqrt{A}z)^{\bar{k}_2-1}\sum_{i=0}^{m}a_{i}^{m}\bar{k}_2\sigma(x^2+Ay^2)^{m-i}(y^2+z^2)^ix\\
&-(x+\sqrt{A}z)^{\bar{k}_2}\sum_{i=0}^{m}a_{i}^{m}2C(m-i)(x^2+Ay^2)^{m-i-1}(y^2+z^2)^ixz\\
&+(x+\sqrt{A}z)^{\bar{k}_2-1}\sum_{i=0}^{m}a_{i}^{m}\bar{k}_0(\sqrt{A}-C)(x^2+Ay^2)^{m-i}(y^2+z^2)^iz\\
&+(x+\sqrt{A}z)^{\bar{k}_2}\sum_{i=0}^{m}a_{i}^{m}2A(m-i)(1-\sigma)(x^2+Ay^2)^{m-i-1}(y^2+z^2)^iy^2.
\end{align*}
Using the transformations (\ref{D9}),
the above equation becomes
\begin{align}\label{DD7}
&\sqrt{u-Aw^2}\sqrt{v-w^2}\frac{d\bar{f}_{2m+\bar{k}_2-1}}{dw}-\bar{k}_2w\bar{f}_{2m+\bar{k}_2-1}\\\nonumber
&=(\sqrt{u-Aw^2}-\sqrt{A}\sqrt{v-w^2})^{\bar{k}_2}\sum_{i=0}^{m}a_{i}^{m}[2\sigma(m-i)+k_{0}+2i]u^{m-i}v^i\\\nonumber
&+(\sqrt{u-Aw^2}-\sqrt{A}\sqrt{v-w^2})^{\bar{k}_2-1}\sum_{i=0}^{m}a_{i}^{m}\bar{k}_2\sigma u^{m-i}v^i\sqrt{u-Aw^2}\\\nonumber
&+(\sqrt{u-Aw^2}-\sqrt{A}\sqrt{v-w^2})^{\bar{k}_2}\sum_{i=0}^{m}a_{i}^{m}2C(m-i)u^{m-i-1}v^i\sqrt{u-Aw^2}\sqrt{v-w^2}\\\nonumber
&+(\sqrt{u-Aw^2}-\sqrt{A}\sqrt{v-w^2})^{\bar{k}_2-1}\sum_{i=0}^{m}a_{i}^{m}\bar{k}_0(C-\sqrt{A})u^{m-i}v^i\sqrt{v-w^2}\\\nonumber
&+(\sqrt{u-Aw^2}-\sqrt{A}\sqrt{v-w^2})^{\bar{k}_2}\sum_{i=0}^{m}a_{i}^{m}2A(m-i)(1-\sigma)u^{m-i-1}v^iw^2,
\end{align}
which is a non-homogeneous linear ordinary differential equation in $\bar{f}_{2m+\bar{k}_2-1}$,
The corresponding homogeneous equation
\begin{eqnarray*}
\sqrt{u-Aw^2}\sqrt{v-w^2}\frac{d\bar{f}^{*}_{2m+\bar{k}_2-1}}{dw}-\bar{k}_2w\bar{f}^{*}_{2m+\bar{k}_2-1}=0
\end{eqnarray*}
has a general solution
\begin{eqnarray*}
\bar{f}^{*}_{2m+\bar{k}_2-1}=(\sqrt{u-Aw^2}-\sqrt{A}\sqrt{v-w^2})^{\bar{k}_2}\bar{A}^{*}_{2m-1}(u,v),
\end{eqnarray*}
where $\bar{A}^{*}_{2m-1}(u,v)$ is an arbitrary smooth function in $u$ and $v$.
In order to use the method of variation of constants,
we assume that
\begin{eqnarray*}
\bar{f}_{2m+\bar{k}_2-1}=(\sqrt{u-Aw^2}-\sqrt{A}\sqrt{v-w^2})^{\bar{k}_2}\bar{A}_{2m-1}(u,v,w)
\end{eqnarray*}
is a solution of (\ref{DD7}),
then $\bar{A}_{2m-1}(u,v,w)$ satisfies
\begin{align*}
\frac{d\bar{A}_{2m-1}}{dw}&=\sum_{i=0}^{m}a_{i}^{m}[2\sigma(m-i)+b+2i]u^{m-i}v^i\frac{1}{\sqrt{u-Aw^2}\sqrt{v-w^2}}\\
&+\sum_{i=0}^{m}a_{i}^{m}\bar{k}_2\sigma u^{m-i}v^i\frac{1}{(\sqrt{u-Aw^2}-\sqrt{A}\sqrt{v-w^2})\sqrt{v-w^2}}\\
&+\sum_{i=0}^{m}a_{i}^{m}2C(m-i)u^{m-i-1}v^i\\
&+\sum_{i=0}^{m}a_{i}^{m}\bar{k}_0(C-\sqrt{A})u^{m-i}v^i\frac{1}{(\sqrt{u-Aw^2}-\sqrt{A}\sqrt{v-w^2})\sqrt{ u-Aw^2}}\\
&+\sum_{i=0}^{m}a_{i}^{m}2A(m-i)(1-\sigma)u^{m-i-1}v^i\frac{w^2}{\sqrt{u-Aw^2}\sqrt{v-w^2}}.
\end{align*}
Some easy computations lead to
\begin{eqnarray}\label{ss}
\int\frac{d w}{(\sqrt{u-Aw^2}-\sqrt{A}\sqrt{v-w^2})\sqrt{u-Aw^2}}
=\frac{w}{u-v}+\int\frac{\sqrt{v-w^2}}{\sqrt{u-Aw^2}}d w.
\end{eqnarray}
Since
\begin{align*}
\int\frac{d w}{\sqrt{u-Aw^2}\sqrt{v-w^2}},~~~~~~~
\int\frac{\sqrt{u-Aw^2}}{\sqrt{v-w^2}}d w
\end{align*}
are elliptic integrals of the first and second kind, by (\ref{D15}) and (\ref{ss}),
in order for $A_{2m-1}(x,y,z)=\bar{A}_{2m-1}(u,v,w)$ to be a homogeneous polynomial
of degree $2m-1$, we must have
\begin{equation}\label{D26}
\begin{cases}
a_{i}^{m}[2\sigma(m-i)+k_0+2i]=0,\\
a_{i}^{m}\bar{k}_{2}\sigma=0,\\
a_{i}^{m}\bar{k}_{2}(C-\sqrt{A})=0,\\
a_{i}^{m}2A(m-i)(1-\sigma)=0, ~~~i=0,1,\cdots,m,
\end{cases}
\end{equation}
which implies $\bar{k}_2=0$, that is to say, $k_2=0$. So the cofactor $K=k_{0}$ is a constant.
Then (\ref{D26}) becomes
\begin{equation*}
\begin{cases}
a_{i}^{m}[2\sigma(m-i)+2i+k_{0}]=0,\\
a_{i}^{m}2A(m-i)(1-\sigma)=0, ~~~i=0,1,\cdots,m.
\end{cases}
\end{equation*}

Now, we split the proof in two cases.
\\
\emph{Case I}. $\sigma=1$, $k_{0}=-2m$. In this case, we have
\begin{align}
&f_{n}=f_{2m}=\sum_{i=0}^{m}a_{i}^{m}(x^2+Ay^2)^{m-i}(y^2+z^2)^{i},\label{ss1}\\
&f_{n-1}=f_{2m-1}=\sum_{i=0}^{m}a_{i}^{m}2C(m-i)(x^2+Ay^2)^{m-i-1}(y^2+z^2)^{i}y.\label{ss2}
\end{align}
Substituting (\ref{ss1}) and (\ref{ss2}) into (\ref{DD4}) with $i=n-2$, we get
\begin{align*}
Ayz\frac{\partial f_{n-2}}{\partial x}&-xz\frac{\partial f_{n-2}}{\partial y}+xy\frac{\partial f_{n-2}}{\partial z}\\
=&-\sum_{i=0}^{m}a_{i}^{m}4C^{2}(m-i)(m-i-1)(x^2+Ay^2)^{m-i-2}(y^2+z^2)^ixyz\\
&-\sum_{i=0}^{m}a_{i}^{m}2Rai(x^2+Ay^2)^{m-i}(y^2+z^2)^{i-1}z\\
&-\sum_{i=0}^{m}a_{i}^{m}4AC(m-i)(x^2+Ay^2)^{m-i-2}(y^2+z^2)^{i}y^{3}.
\end{align*}
By (\ref{D9}),
the above equation becomes
\begin{align*}
\frac{d{\bar{f}}_{n-2}}{dw}=&\sum_{i=0}^{m}a_{i}^{m}4C^{2}(m-i)(m-i-1)u^{m-i-2}v^iw\\
&+\sum_{i=0}^{m}a_{i}^{m}2Raiu^{m-i}v^{i-1}\frac{1}{\sqrt{u-Aw^{2}}}\\
&-\sum_{i=0}^{m}a_{i}^{m}2C(m-i)u^{m-i-1}v^{i}\frac{w}{\sqrt{u-Aw^{2}}\sqrt{v-w^{2}}},
\end{align*}
that is,
\begin{align*}
\bar{f}_{n-2}&=\sum_{i=0}^{m}a_{i}^{m}2C^{2}(m-i)(m-i-1)u^{m-i-2}v^i{w^2}\\
&+\sum_{i=0}^{m}a_{i}^{m}2Raiu^{m-i}v^{i-1}\frac{1}{\sqrt{A}}\arctan\left(\frac{\sqrt{A}w}{\sqrt{u-Aw^2}}\right)\\
&-\sum_{i=0}^{m}a_{i}^{m}2C(m-i)u^{m-i-1}v^{i}\frac{\ln|\sqrt{u-Aw^2}+\sqrt{A(v-w^2)}|}{\sqrt{A}}.
\end{align*}
We must have
\begin{equation*}
\begin{cases}
a_{i}^{m}2i Ra=0,\\
a_{i}^{m}2C(m-i)=0,
 ~~~i=0,1,\cdots,m,
\end{cases}
\end{equation*}
which yields $
a_{i}^{m}=0$ for $i=0,1,\cdots,m.
$
This is a contradiction to $f_n\neq0$.
\\
\\
\emph{Case II}. $\sigma\neq1$, $a_{i}^{m}=0$, $i=0,1,\cdots,m-1$, $k_{0}=-2\sigma m$. In this case, we have
\begin{align*}
f_{n}=f_{2m}=a_{m}^{m}(y^2+z^2)^{m},~~f_{n-1}=f_{2m-1}=0.
\end{align*}
Working in a
similar way to the previous case, we get
\begin{align*}
-\sqrt{u-Aw^{2}}\sqrt{v-w^{2}}\frac{d\bar{f}_{n-2}}{dw}
=-2a_{m}^{m}mRav^{m-1}\sqrt{v-w^{2}},
\end{align*}
and
\begin{align*}
\bar{f}_{n-2}
=2a_{m}^{m}mRav^{m-1}\frac{1}{\sqrt{A}}\arctan\left(\frac{\sqrt{A}w}{\sqrt{u-Aw^2}}\right).
\end{align*}
Therefore, $a_{m}^{m}=0$  and $f_n=0$, which is a contradiction.
\end{proof}

\begin{prop}\label{T3}
System (\ref{I1}) has no exponential factors.
\end{prop}
\begin{proof}
Let $E = \exp(g/h)$ be an
exponential factor of system (\ref{I1}) with a cofactor
$$
L=l_0+l_1x+l_2y+l_3z,~~l_i\in\mathbb{R},~~i=0,1,2,3,
$$
 where $g,h\in\mathbb{R}[x,y,z]$ with
$f,g$ being prime. From Proposition \ref{PP2}, \ref{T1} and \ref{T2},  $E = \exp(g)$ with $g=g(x,y,z)\in \mathbb{R}[x,y,z]/{\mathbb{R}}$.
By definitions,  $g$ satisfies
\begin{align}\label{E1}
(Ayz+Cz-\sigma x)\frac{\partial g}{\partial x}+(-xz+Ra-y)\frac{\partial g}{\partial y}&+(-z+xy)\frac{\partial g}{\partial z}=l_{0}+l_{1}x+l_{2}y+l_{3}z.
\end{align}

If  $g$ is a polynomial of degree $n\geq3$.
We write $g$ as $g=\sum_{j=0}^{n}g_{j}(x,y,z)$,
where each $g_{j}$ is a homogeneous polynomial of degree $j$ and $g_n\neq 0$.
Observing the right hand side of (\ref{E1}) has degree at most one,
we compute the terms of degree $n+1$ in (\ref{E1}) and get
\begin{equation*}\label{E2}
Ayz\frac{\partial g_{n}}{\partial x}-xz\frac{\partial g_{n}}{\partial y}+xy\frac{\partial g_{n}}{\partial z}=0,
\end{equation*}
which is (\ref{po}) replacing $f_n$ by $g_n$. Then the arguments used in the proof of Proposition \ref{T1} imply
$g_{n}=L_{n}(y^{2}+z^{2})^{m}$ with $L_{n}\in\mathbb{R}$.
Now computing the terms of degree $n$ in (\ref{E1}) leads to
\begin{align*}\label{E3}
Ayz\frac{\partial g_{n-1}}{\partial x}-xz\frac{\partial g_{n-1}}{\partial y}+xy\frac{\partial g_{n-1}}{\partial z}&=(\sigma x-Cz)\frac{\partial g_{n}}{\partial x}+y\frac{\partial g_{n}}{\partial y}+z\frac{\partial g_{n}}{\partial z},
\end{align*}
which is  (\ref{assd}) with $f_n$ replaced by $g_n$ and $f_{n-1}$ replaced by $g_{n-1}$. Again, the arguments used in the
proof of Proposition \ref{T1} imply that $g_n=0$, which is a contradiction.

Hence, $g$ is a polynomial of degree at most two satisfying (\ref{E1}).
So we can write $g$ as
\begin{equation}\label{E4}
g=b_{0}+b_{1}x+b_{2}y+b_{3}z+b_{11}x^{2}+b_{22}y^{2}+b_{33}z^{2}+b_{12}xy+b_{13}xz+b_{23}yz.
\end{equation}
Then substituting (\ref{E4}) into (\ref{E1}) yields a system of  algebraic equations
 $g_i=0$  for $i=1,\cdots,15$, where
\begin{align*}
&g_1=Ra b_{2}-l_0,~g_2=-\sigma b_{1}+Ra b_{12}-l_1,\\
&g_3=2Ra b_{22}- b_{2}-l_2,~g_4=Cb_{1}+Ra b_{23}-b_3-l_3,\\
&g_5=-\sigma b_{12}-b_{12}+b_{3},~g_6=-\sigma b_{13}-b_{2}-b_{13},\\
&g_7=A b_{1}+Cb_{12}-2b_{23},~g_8=-2\sigma b_{11},\\
&g_9=-2b_{22},~g_{10}=2Cb_{11}+Cb_{13}-2b_{33},~g_{11}=-b_{22}+b_{33},\\
&g_{12}=b_{13},
~g_{13}=b_{12},
~g_{14}=b_{23},
~g_{15}=2Ab_{11}+Ab_{13}.
\end{align*}
From the above equations,
we obtain that $l_{0}=l_{1}=l_{2}=l_{3}=0$, $b_{1}=b_{2}=b_{3}=0$, $b_{11}=b_{22}=b_{33}=b_{12}=b_{13}=b_{23}=0$,
which implies $g$ is a constant.
This completes the proof.
\end{proof}

 By Proposition \ref{PP1} and Proposition \ref{T1}-\ref{T3}, we obtain the main result in this section.
\begin{thm}\label{th}
The following statements hold for the Glukhovsky-Dolzhansky system (\ref{I1}).\\
(a) It has no polynomial first integrals.\\
(b) It has no Darboux polynomials with non-zero cofactors.\\
(c) It has no exponential factors.\\
(d) It admits no Darboux first integrals.
\end{thm}

\begin{remark}
From a physical point of view,
the parameters in the Glukhovsky-Dolzhansky system (\ref{I1}) are positive due to their physical meaning.
Whereas,
from a mathematical point of view,
one may consider non-positive parameters and then some integrable cases of system (\ref{I1}) are founded:
\begin{itemize}
\item $Ra=0$ and $\sigma=1$, system (\ref{I1}) has a Darboux polynomial $f=y^2+z^2$ with a constant cofactor $k=-2$.
\item $C=0$ and $\sigma=1$, system (\ref{I1}) has a Darboux polynomial $f=x^2-Az^2$ with a constant cofactor $k=-2$.
\item $Ra=C=0$ and $\sigma=1$, system (\ref{I1}) has a rational first integral $\Phi=(y^2+z^2)/(x^2-Az^2)$.
\end{itemize}
But we do not have a clear physical understanding of the above integrable results for the Glukhovsky-Dolzhansky system (\ref{I1}).

\end{remark}

\begin{remark}
Another tool to study the non-integrability of non-Hamiltonian systems is  the differential Galois theory \cite{ga1,ga2,ga3}.
Observing system (\ref{I1}) has a straight line solution $(x(t),y(t),z(t))=(0,Ra-e^{-t},0)$,
we can analyze the differential
Galois group of the normal variational equations along this solution,
and show that system (\ref{I1}) is not rationally integrable in Bogoyavlenskij sense for almost all parameter values, see \cite{gaa1,gaa2,gaa3,newgaa3}
for more details.
\end{remark}

\section{Stability and bifurcations of the GD model}

In this section,
we first point out some basic dynamical properties of the GD model.
The divergence of the vector field $\mathcal{X}$ corresponding to
the GD model is given
$$
\textit{div}~\mathcal{X}=\frac{\partial \dot x}{\partial x}+\frac{\partial \dot y}{\partial y}+\frac{\partial \dot z}{\partial z}=-\sigma-2<0,
$$
which implies that the GD model is dissipative.
The GD model admits the symmetry $(x,y,z)\rightarrow (-x,y,-z)$,
that is, $(x(t),y(t),z(t))$ is a solution of the GD model if and only if $(-x(t),y(t),-z(t))$ is also a solution of that.
In addition,
we show that the GD model possesses a global attractor as follows.
\begin{thm}\label{global}
There exits a positive number $M$ such that the ellipsoid
$$
S_M=\left\{(x,y,z)\bigg|\frac{\sigma x^2}{2}+\frac{(A\sigma +C^2)y^2}{2}+\frac{C^2z^2}{2}\leq M\right\}
$$
is a global attractor of the GD model, that is, any solution $(x(t),y(t),z(t))$ with
the initial condition $(x(0),y(0),z(0))$ outside $V_m$ will enter $V_m$ as $t\rightarrow +\infty$.
Thus,
the GD model is ultimately bounded.
\end{thm}
\begin{proof}
Let
$$
V(x,y,z)=\frac{\sigma x^2}{2}+\frac{(A\sigma +C^2)y^2}{2}+\frac{C^2z^2}{2}.
$$
Then the rate of change of $V$ along a solution of the GD model reads
\begin{align*}
\frac{dV}{dt}=&\dot x\frac{\partial V}{\partial x}+\dot y\frac{\partial V}{\partial y}+\dot z\frac{\partial V}{\partial z}\\
=&\sigma C xz-\sigma^2x^2-C^2z^2-(A\sigma+C^2)(y-\frac{R_a}{2})^2+\frac{R_a^2(A\sigma+C^2)}{4}\\
\leq&\sigma C(\frac{\sigma x^2}{2C}+\frac{Cy^2}{2\sigma})-\sigma^2x^2-C^2z^2-(A\sigma+C^2)(y-\frac{R_a}{2})^2+\frac{Ra^2(A\sigma+C^2)}{4}\\
=&-\frac{\sigma^2x^2}{2}-\frac{C^2z^2}{2}-(A\sigma+C^2)(y-\frac{R_a}{2})^2+\frac{R_a^2(A\sigma+C^2)}{4}.
\end{align*}
Define a ellipsoid
$$
S=\left\{(x,y,z)\bigg|\frac{\sigma^2x^2}{2}+\frac{C^2z^2}{2}+(A\sigma+C^2)(y-\frac{R_a}{2})^2\leq\frac{R_a^2(A\sigma+C^2)}{4}\right\}
$$
Then, we can select a sufficiently large number $M>0$ such that $S\subset S_M$.
Hence, for any point $(x,y,z)$ outside $S_M$ we find $\frac{dV}{dt}|_(x,y,z)<0$. This completes the proof.
\end{proof}

The next result gives the existence of equilibrium points.
\begin{prop}
There is a threshold $R_a^*=R_a^*(\sigma,A,C)$ such that

(i) if $0<R_a\leq R_a^*$, the GD model has only one  equilibrium point $E_0=(0,R_a,0)$;

(ii) if $R_a>R_a^*$, the GD model has three and only three equilibrium points $E_0$,
$$
E_{\pm}=(\pm\sqrt{\eta},\frac{R_a}{\eta+1},\pm\frac{R_a\sqrt{\eta}}{\eta+1}),~~\textit{wherer}~~\eta=\frac{CR_a-2\sigma+R_a\sqrt{C^2+4AC}}{2\sigma}.
$$
\end{prop}
\begin{proof}
Applying $\dot x=\dot y=\dot z=0$,
we obtain the equilibrium points satisfying
$$
Ayz+Cz-\sigma x=0,~-xz+R_a-y=0,~-z+xy=0.
$$
Eliminating the variables $y,~z$ leads to
$$
y=\frac{R_a}{x^2+1},~~~~~~~~z=\frac{xR_a}{x^2+1},
$$
and
\begin{align}\label{xeq}
\sigma x^5+(2\sigma-CR_a)x^3+(\sigma-CR_a-AR_a^2)x=0.
\end{align}
It is easy to check that if $\sigma-CR_a-AR_a^2\geq0$,  (\ref{xeq}) has only one real root $x=0$; if  $\sigma-CR_a-AR_a^2<0$, (\ref{xeq}) has three real roots $x=0,\pm\sqrt{\eta}$. Define the threshold value
$$
R_a^*=\frac{\sqrt{C^2+4A\sigma}-C}{2A}.
$$
We can complete the proof easily.
\end{proof}

Using the standard Routh-Hurwitz criterion, we also obtain the local stability of each  equilibrium points.
Let
\begin{align*}
L(A,\sigma,C)=&-AC\sigma^3+3AC\sigma^2+C^3\sigma-6AC\sigma-2C^3\\
&+\sqrt{4A\sigma+C^2}(-A\sigma^3+A\sigma^2+C^2\sigma-2A\sigma-2C^2),\\
R(A,\sigma,C)=&C^2\sigma^3+8A\sigma^3+4C^2\sigma^2+\sqrt{4A\sigma+C^2}(C\sigma^3+4C\sigma^2).
\end{align*}
We also define a threshold $R_a^{**}=L(A,\sigma,C)/R(A,\sigma,C)$.
\begin{prop}\label{sta}
 The following statements hold.

(i) $E_0$ is  asymptotically stable if and only if $R_a<R_a^*$.

(ii) Supposing $R_a>R_a^*$,  $E_{\pm}$ is  asymptotically stable if and only if either $L(A,\sigma,C)\leq 0$ or
$L(A,\sigma,C)>0$, $R_a<R_a^{**}$.

\end{prop}
\begin{proof}
(i) Some direct computations yield
the  characteristic polynomial at $E_0$ is
\begin{align}
P_0(\lambda)&=\lambda^3+p_{11}\lambda^2+p_{12}\lambda+p_{13}\notag\\
   &=\lambda^3+(\sigma+2)\lambda^2+(-AR_a^2-CR_a+2\sigma+1)\lambda-AR_a^2-CR_a+\sigma. \label{e0ch}
\end{align}
Clearly, we have
\begin{align}\label{con}
p_{11}>2,~~p_{12}>p_{13}+1.
\end{align}
Combing with the Routh-Hurwitz criterion (see Proposition A1 in Appendix),
we know that the real parts of all the roots $\lambda$ are negative if and only if $p_{13}>0$,
i.e., $R_a<R_a^*$.

(ii) Similarly, the  characteristic polynomial at $E_{\pm}$ is
\begin{align}
P_{\pm}(\lambda)&=\lambda^3+p_{21}\lambda^2+p_{22}\lambda+p_{23}\notag\\
   &=\lambda^3+(\sigma+2)\lambda^2+\frac{R_a^2(M_1R_a+M_2)}{2\sigma^3(\eta+1)^2}\lambda+\frac{R_a^2(M_3R_a-M_4)}{\sigma^2(\eta+1)^2},\label{cha2}
\end{align}
where the positive numbers $M_i$, $i=1,2,3,4$ are given by
\begin{align*}
M_1=&AC\sigma^2+A\sigma^2\sqrt{4A\sigma+C^2}+3AC\sigma+A\sigma\sqrt{4A\sigma+C^2}+C^3+C^2\sqrt{4A\sigma+C^2},\\
M_2=&C^2\sigma^2+C\sigma^2\sqrt{4A\sigma+C^2},\\
M_3=&4AC\sigma+2A\sigma\sqrt{4A\sigma+C^2}+C^3+C^2\sqrt{4A\sigma+C^2},\\
M_4=&4A\sigma^2+C^2\sigma+C\sigma\sqrt{4A\sigma+C^2}.
\end{align*}
Observing $R_a^*=M_4/M_3$, we see that
\begin{align}\label{cond2}
p_{21}>2,~p_{22}>0,~
p_{23}>\frac{R_a^2(M_3R_a^*-M_4)}{\sigma^2(\eta+1)^2}=0.
\end{align}
Hence, by the Routh-Hurwitz criterion (see Proposition A1 in Appendix),
we know that the real parts of all the roots $\lambda$ of (\ref{cha2}) are negative if and only if
$
p_{21}p_{22}>p_{23},
$
or equivalently
$
L(A,\sigma,C)R_a<R(A,\sigma,C).
$
\end{proof}
\begin{remark}
One can also easily prove that $E_0$ is global asymptotically stable if $R<R^*$.
Indeed, supposing $R<R^*$ (or equivalently $\sigma>R_a(AR_a+C)$),
consider the Lyapunov function
$$
V_{global}(x,y,z)=\frac{x^2}{2}+\frac{(2\sigma-CR_a)}{2R_a^2}(y-R_a)^2+\frac{2\sigma-R_a(AR_a+C)}{2R_a^2}z^2.
$$
We have
$$
\frac{dV}{dt}=-\sigma\left(x-\frac{z}{R_a}\right)^2-\frac{\sigma-R_a(AR_a+C)}{R_a^2}z^2-\frac{2\sigma-CR_a}{R_a^2}(y-R_a)^2<0
$$
for any $(x,y,z)\neq E_0$.
It follows that $E_0$ is global asymptotically stable under the condition
$R<R^*$.
\end{remark}

\begin{figure}[htpp]
\centering
\label{xx22}
\includegraphics[width=0.45\textwidth]{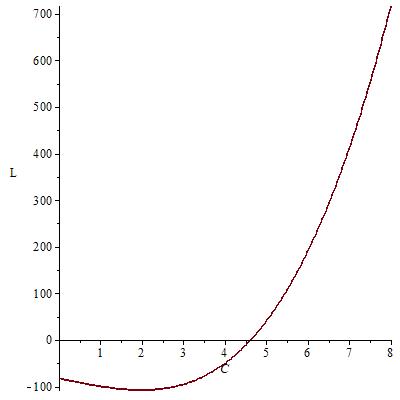}
\caption{The graph of the function $L=L(1,3,C)$.}
\label{Fig.lable}
\end{figure}

\begin{remark}
We remark that the parameter space $(A,\sigma,C)\in (\mathbb{R}^+)^3$ satisfying $L(A,\sigma,C)>0$ is not empty.
For example, fixed $A=1,~\sigma=3$ the graph of the function $L=L(1,3,C)$ is shown in Figure \ref{xx22}.
\end{remark}

Now we turn to study the bifurcations of the stead states $E_0$ and $E_{\pm}$.
A necessary condition for the occurrence of bifurcations is either a simple real eigenvalue approaches
zero or a pair of simple complex eigenvalues reaches the imaginary axis for some values of the parameters.
Recall that all co-dimension one or two bifurcations of equilibria are classified as follows.
\\
\\
 Co-dimension one bifurcations
\begin{itemize}
\item {\bf Fold bifurcation} corresponding to  the equilibria has
a simple zero eigenvalue.

\item  {\bf Hopf bifurcation of co-dimension one} corresponding to the equilibria has a simple pair of purely imaginary eigenvalues.
\end{itemize}

 Co-dimension two bifurcations
\begin{itemize}
\item {\bf Cusp bifurcation} corresponding to the equilibria has a simple zero eigenvalue and
the quadratic coefficient vanishes.

\item {\bf Hopf bifurcation of co-dimension two} (also called \emph{generalized Hopf bifurcation}) corresponding to the equilibria has a simple pair of purely imaginary eigenvalues and the first Lyapunov coefficient vanishes.

\item {\bf Bogdanov-Takens bifurcation} corresponding to the equilibria has a zero eigenvalue of (algebraic) multiplicity two.

\item {\bf Fold-Hopf bifurcation} corresponding to the equilibria has a simple zero eigenvalue and a
 simple pair of purely imaginary eigenvalues.

\item {\bf Double Hopf bifurcation} corresponding to the equilibria has two pairs of purely imaginary eigenvalues.
\end{itemize}

Clearly, the GD model cannot undergo double Hopf bifurcations since it has only three eigenvalues.
We first investigate bifurcations of the GD model at $E_0$.
The characteristic polynomial at $E_0$ is given in (\ref{e0ch}).
By Proposition \ref{ap1} and (\ref{con}),
we obtain that the possible bifurcations of co-dimension one or two are the fold bifurcation or
cusp bifurcation.

Next, we use the center manifold reduction method to determine the type of bifurcations at $E_0$.
We translate $E_0$ to the origin by $(x_1,y_1,z_1)=(x,y-Ra,z)$ and rewrite the GD model into
\begin{align}\label{gd2}
\dot X=F(X)=JX+\frac{1}{2}B(X,X), ~~~~~~~X=(x_1,y_1,z_1),
\end{align}
where
$$
J=\left(
    \begin{array}{ccc}
      -\sigma & 0 & AR_a+C \\
      0 & -1 & 0 \\
      R_a & 0 & -1 \\
    \end{array}
  \right),~~B(x_1,y_1,z_1,x_2,y_2,z_2)=\left(
                                         \begin{array}{c}
                                         A(y_1z_2+y_2z_1) \\
                                           -x_1z_2-x_2z_1\\
                                           x_1y_2+x_2y_1 \\
                                         \end{array}
                                       \right).
$$
When $\sigma=AR_a^2+CR_a$ (or equivalently $R_a=R_a^*$), the Jacobian matrix $J$ has a simple zero  eigenvalue.
Define the vectors
$$
q=\left(
    \begin{array}{c}
      \frac{1}{R_a} \\
      0 \\
      1 \\
    \end{array}
  \right),~~
  p=\left(
      \begin{array}{c}
        \frac{R_a}{(AR_a^2+CR_a+1} \\
        0 \\
       \frac{R_a(AR_a+C)}{(AR_a^2+CR_a+1}
      \end{array}
    \right)
$$
such that $Jq=0,~J^Tp=0$ and $\langle p,q\rangle=1$.
The one-dimensional center manifold can
be parameterized by $w$ through the immersion
$$
X=H(w)=wq+\frac{1}{2}h_2w^2+\frac{1}{6}h_3w^3+O(w^4),~~h_2,h_3\in\mathbb{R}^3,
$$
and  system (\ref{gd2}) can be
reduced to
$$
\dot w=G(w)=bw^2+cw^3+O(w^4),~~b,c\in\mathbb{R}.
$$
The unknown quantities $h_2,~h_3,~b,~c$ are determined by the \emph{homological equation} \cite{bifur}
$$
H_w(w)G(w) =F(H(w)).
$$
More precisely, in our case we have
$$
b=\frac{1}{2}\langle p,B(q,q)\rangle=0,~c=\frac{1}{6}\langle p,3B(p,3B(q,h_2))\rangle=-\frac{2AR_a^*+C}{R_a^*((AR_a^*)^2+CR_a^*+1))}<0,
$$
where $h_2=(0,-2/R_a^*,0)^T$ is the  solution of
$$
Jh_2+B(q,q)=0,~~\langle p,h_2\rangle=0.
$$
Then generically the GD model restricted to the center manifold is locally topologically equivalent to the normal form
$$
\dot w=(\beta_1+\beta_2w)w-\frac{2AR_a^*+C}{R_a^*((AR_a^*)^2+CR_a^*+1))}w^3,
$$
with $\beta_1,\beta_2$ being unfolding parameters.

 To summary, we get the result on the cusp bifurcation of the GD model at $E_0$. Recall that the cusp bifurcation implies
the presence of a hysteresis phenomenon \cite{bifur2}.
\begin{thm}
Assume $R_a=R_a^*$. Then the GD model undergoes a cusp bifurcation  at $E_0$.
\end{thm}

Proceeding as above, we investigate the bifurcations of the GD model at $E_{\pm}$.
As mentioned above, the GD model is invariant under the transformation
$$
(x, y, z) \rightarrow (-x,y,-z).
$$
This implies that dynamics in the neighbourhood of a
point $(x, y, z)$ looks the same as that of a point $(-x,y,-z)$.
Hence, it is sufficient to study the bifurcations occurring at $E_{+}$ which is also true for $E_{-}$.
Using Proposition A2 in Appendix, (\ref{cha2}) and (\ref{cond2}),
we see that the possible bifurcations of co-dimension one or two at $E_{+}$ are Hopf bifurcations.

If the physical parameters $A,~\sigma,~C$ satisfy the nonlinear inequality $L(A,\sigma,C)\leq 0$.
By Proposition \ref{sta},
$E_+$ is asymptotically stable and cannot undergo bifurcation.
So we assume $L(A,\sigma,C)$ is positive.
The next result gives the existence of Hopf bifurcation at $E_{+}$.

\begin{thm}
Assume $L(A,\sigma,C)>0$ and $R_a=R_a^{**}$. The GD model undergoes a Hopf bifurcation at $E_{+}$.
\end{thm}
\begin{proof}
For $L(A,\sigma,C)>0$ and $R_a=R_a^{**}$,
it follows that (\ref{cha2}) has one negative real root $\lambda_1=-\sigma-2$ and a pair of conjugate purely imaginary
roots
$$
\lambda_{2,3}=\pm w_0i,~w_0=R_a^{**}\sqrt{\frac{M_1R^{**}_a+M_2}{2\sigma^3(\eta+1)^2}}.
$$
Substituting $\lambda_{2,3}=\lambda_{2,3}(R_a)$ into (\ref{cha2}) and taking the derivative with respect to $R_a$,
we have
$$
\frac{d\lambda_{2,3}}{dR_a}=-\left(\frac{dp_{21}}{dR_a}\lambda_{2,3}^2+{\frac{dp_{22}}{dR_a}\lambda_{2,3}+\frac{dp_{23}}{dR_a}}\right)\frac{1}{3\lambda_{2,3}^2+2p_1\lambda_{2,3}+p_2}.
$$
Combing with the real part $\Re(\lambda_{2,3})=0$ at $R_a=R_a^{**}$,
we get that
\begin{align*}
\frac{d\Re(\lambda_{2,3})}{dR_a}\bigg|_{R_a=R_a^{**}}&=\left(p_{21}\frac{dp_{22}}{dR_a}-p_{22}\frac{dp_{21}}{dR_a}+\frac{dp_{23}}{dR_a}\right)\frac{1}{2p_{22}+2p_{21}^2}\bigg|_{R_a=R_a^{**}}.
\end{align*}
Some tedious calculations yield the explicit expression of the numerator
$$
p_{21}\frac{dp_{22}}{dR_a}-p_{22}\frac{dp_{21}}{dR_a}+\frac{dp_{23}}{dR_a}=\frac{
2A\sigma^3+14A\sigma^2+3C^2\sigma+3C\sqrt{4A\sigma+C^2}\sigma+4A\sigma+2C^2+2C\sqrt{4A\sigma+C^2}}{\sigma(C+\sqrt{4A\sigma+C^2})}.
$$
Hence
$$
\frac{d\Re(\lambda_{2,3})}{dR_a}\bigg|_{R_a=R_a^{**}}>0
$$
is positive,
and the transversal condition holds.
All conditions for a Hopf bifurcation to occur are met \cite{bifur2},
and a Hopf bifurcation happens at $E_{+}$.
\end{proof}
Let us mention that to determine the stability of periodic solution from the Hopf bifurcations,
one can easily apply the center manifold reduction \cite{bifur} to calculate Lyapunov coefficients.
We do not give Lyapunov coefficients here since the expressions of Lyapunov coefficients are too long.

\begin{figure}[htpp]
\centering
\subfigure[$t-x$ space]{
\includegraphics[width=0.4\textwidth]{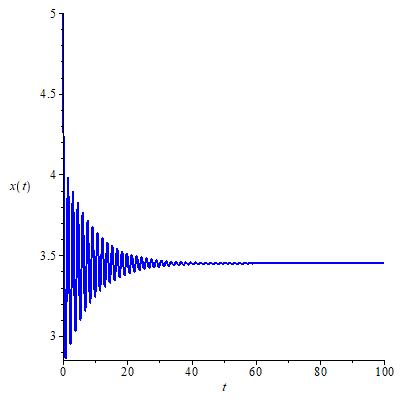}}
\subfigure[$t-y$ space]{
\includegraphics[width=0.4\textwidth]{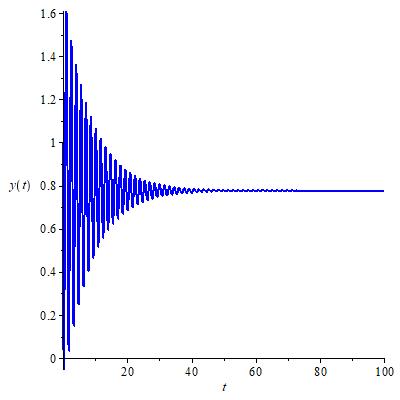}}
\\
\subfigure[$t-z$ space]{
\includegraphics[width=0.4\textwidth]{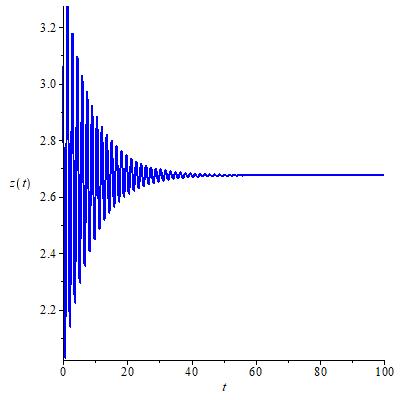}}
\subfigure[$x-y-z$ space]{
\includegraphics[width=0.45\textwidth]{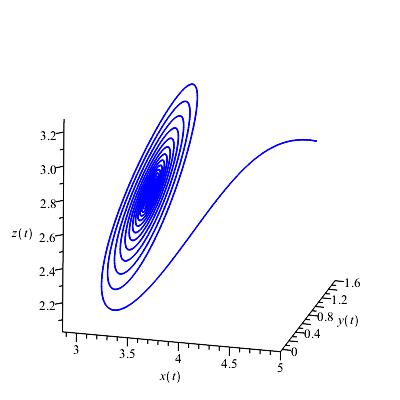}}
\caption{The equilibrium $E_{+}$ is asymptotically stable for the GD model with $R_a=10$.}
\label{Figchaos1}
\end{figure}

\begin{figure}[htpp]
\centering
\subfigure[$t-x$ space]{
\includegraphics[width=0.40\textwidth]{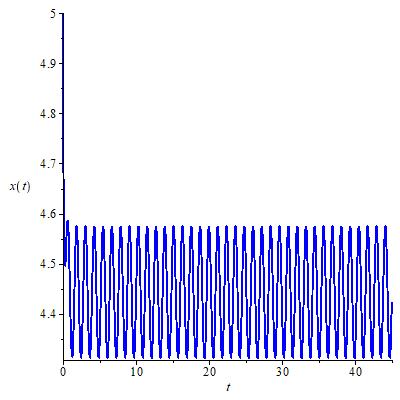}}
\subfigure[$t-y$ space]{
\includegraphics[width=0.4\textwidth]{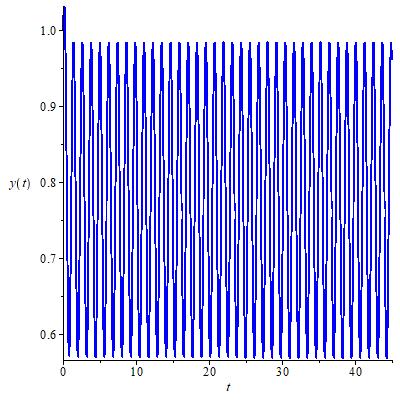}}
\\
\subfigure[$t-z$ space]{
\includegraphics[width=0.40\textwidth]{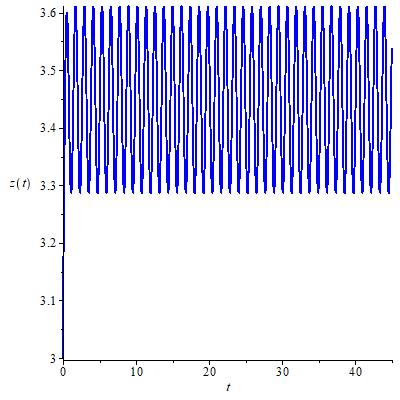}}
\subfigure[$x-y-z$ space]{
\includegraphics[width=0.45\textwidth]{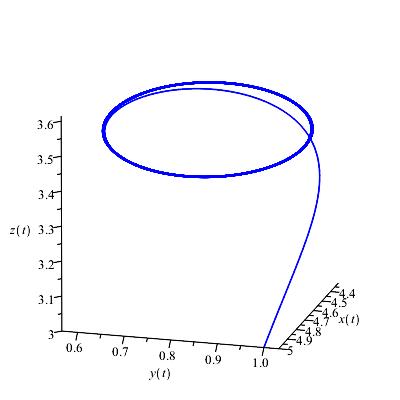}}
\caption{Bifurcating periodic solution of the GD model at $E_{+}$
with $R_a=16.1$}
\label{Figchaos2}
\end{figure}

\begin{figure}[htpp]
\centering
\subfigure[Projections into the plane $(x,y)$]{
\includegraphics[width=0.40\textwidth]{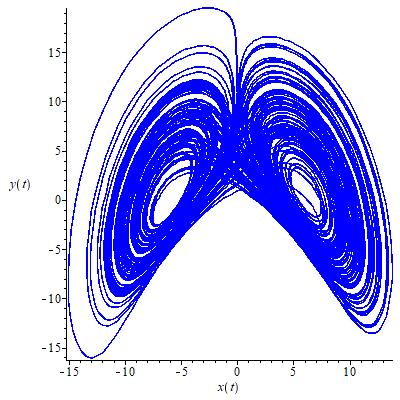}}
\subfigure[Projections into the plane $(x,z)$]{
\includegraphics[width=0.4\textwidth]{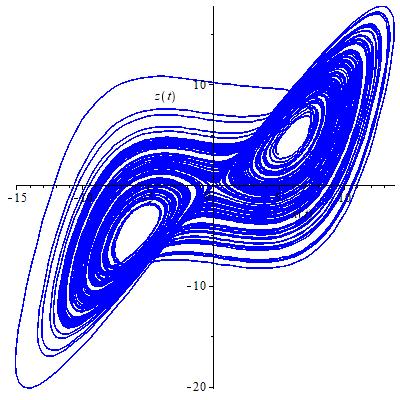}}
\\
\subfigure[Projections into the plane $(y,z)$]{
\includegraphics[width=0.40\textwidth]{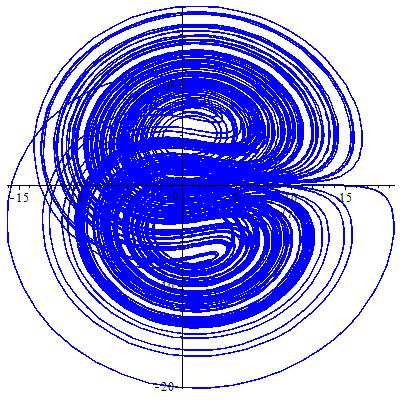}}
\subfigure[$x-y-z$ space]{
\includegraphics[width=0.5\textwidth]{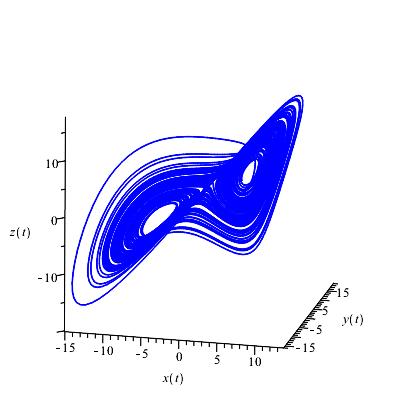}}
\caption{Chaotic attractor of the GD model with $R_a=30$.}
\label{Figchaos3}
\end{figure}

At the end of this section,
we shall give some numerical simulations to verify the above theoretical results.
We choose the parameter values $(A,\sigma,C)=(0.2,4,5)$ such that $L(0.2,4,5)>0$.
Then the threshold values $R_a^*\approx 0.77591805$ and $R_a^{**}\approx 16.10662151$.
Set the initial values
$(x(0),y(0),z(0))=(5,1,3)$.
For small $R_a$, e.g. $R_a=10$,
it follows from Proposition \ref{sta} that $E_+$ is asymptotically stable as shown in Figure \ref{Figchaos1}.
As $R_a$ passes though $R_a^{**}$,
a periodic solution bifurcates from the Hopf bifurcation as shown in Figure \ref{Figchaos2}.
When $R_a$ is much larger than $R_a^{**}$,
the GD model may display some chaotic behaviors,
see Figure \ref{Figchaos3}.

\section{Dynamics at infinity of the GD model}

To understand the final evolution of the orbit $\{x(t),y(t),z(t)\}$ of the GD model,
we need to characterize its global dynamical behaviors, i.e.,
the flow of the GD model at infinity,
which gives all the information about infinite equilibrium points.
To this end,
in this section we will analyze the Poincar\'{e} compactification of the GD model in the local charts $U_{i}$ and $V_{i}$ ($i = 1,2,3$).
We refer to \cite{MR0252788}, \cite{MR998352}, \cite{MR3474921}, \cite{Liu2012Dynamics} for the detailed theory of Poincar\'{e} compactification.
\begin{figure}
  \centering
  \includegraphics[width=230pt]{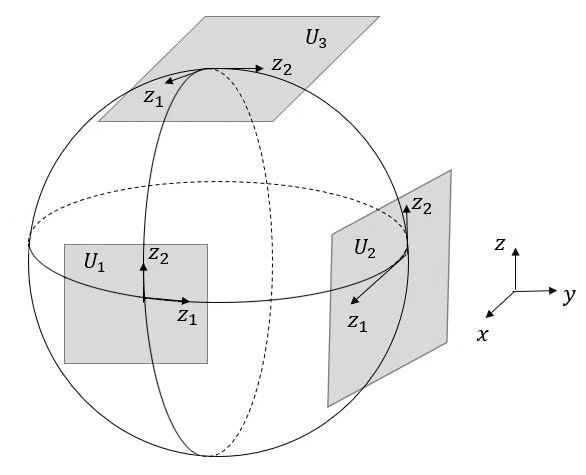}
  \caption{Orientation of the local charts $U_{i}$ and $V_{i}$ ($i = 1,2,3$) in the positive endpoints of the $x$, $y$ and $z$ axis.}
  \label{figure100}
\end{figure}

Let $S^{3}=\{r=(r_{1},r_{2},r_{3},r_{4})\in \mathbb{R}^{4}|\|r\|=1\}$ be a Poincar\'{e} unit sphere \cite{MR0252788}.
We divide the above sphere into three parts,
that is the northern hemisphere $S_{+}=\{r\in S^{3},r_{4}>0\}$,
the southern hemisphere $S_{-}=\{r\in S^{3},r_{4}<0\}$ and the equator $S^{1}=\{r\in S^{3},r_{4}=0\}$.
Denote the tangent hyperplane at the point $(\pm1, 0, 0, 0)$, $(0, \pm1, 0, 0)$, $(0, 0, \pm1, 0)$ and $(0, 0, 0, \pm1)$
by the local chart $U_{i}$ and $V_{i}$ ($i = 1,2,3,4$) where $U_{i}=\{r\in S^{3},r_{i}>0\}$ and $V_{i}=\{r\in S^{3},r_{i}<0\}$.
Define the central projections $f^{+}:\mathbb{R}^{3}\rightarrow S^{3}$ and $f^{-}:\mathbb{R}^{3}\rightarrow S^{3}$ by
\begin{equation}\nonumber
f^{\pm}(x,y,z)=\pm(\frac{x}{\Delta}, \frac{y}{\Delta}, \frac{z}{\Delta}, \frac{1}{\Delta}),~\textrm{where}~\Delta=\sqrt{1+x^{2}+y^{2}+z^{2}},
\end{equation}
and define $\varphi_{k}: U_{k}\rightarrow \mathbb{R}^{3}$, $\phi_{k}: V_{k}\rightarrow \mathbb{R}^{3}$ by $\varphi_{k}=- \phi_{k}=(\frac{r_{l}}{r_{m}}, \frac{r_{m}}{r_{k}}, \frac{r_{n}}{r_{k}})$ for $k=1,2,3,4$ with $1\leq l, m, n\leq4$ and $l,m,n\neq k$.
We only consider the local charts  $U_{i}$ and $V_{i}$ for
$i = 1,2,3$ to get the dynamics at infinity as shown in Figure \ref{figure100}.
\begin{figure}
  \centering
  \includegraphics[width=230pt]{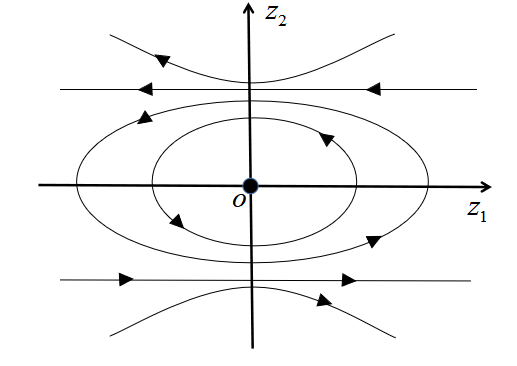}
  \caption{Dynamics of the GD model on the Poincar\'{e} sphere at infinity in the local chart $U_1$}
  \label{u1}
\end{figure}
\emph{In the local charts $U_1$}
Doing the change of variables $(x,y,z)=(z_3^{-1},z_1z_3^{-1},z_2z_3^{-1})$ and rescaling the time ${\rm d}\tau=z_3^{-1}{\rm d}t$,
the Poincar\'{e} compactification $p(X)$ of the GD model in the local chart $U_{1}$ is given by
\begin{equation}\label{nu1}
\begin{cases}
\frac{\textrm{d}z_{1}}{\textrm{d}{\tau}}=-Az_1^2z_2-Cz_1z_2z_3+R_az_3^2+\sigma z_1z_3-z_1z_3-z_2\\
\frac{\textrm{d}z_{2}}{\textrm{d}{\tau}}=-Az_1z_2^2-Cz_2^2z_3+\sigma z_2z_3-z_2z_3+z_1\\
\frac{\textrm{d}z_{3}}{\textrm{d}{\tau}}=-z_3(Az_1z_2+Cz_2z_3-\sigma z_3).
\end{cases}
\end{equation}
Recall that the $z_{1}z_{2}$-plane is invariant under the flow of system (\ref{nu1}),
which completely describes the dynamics on the sphere at infinity.
For $z_{3}=0$, system (\ref{nu1}) reduces to
\begin{equation}\label{nu11}
\frac{\textrm{d}z_{1}}{\textrm{d}{\tau}}=-Az_1^2z_2-z_2,~~~~
\frac{\textrm{d}z_{2}}{\textrm{d}{\tau}}=-Az_1z_2^2+z_1.
\end{equation}
One can check that system (\ref{nu11}) is integrable with a polynomial first integral
$$
\Phi(z_1,z_2)=\frac{Az_2^2-1}{Az_1^2+1}.
$$
Using this first integral,
we have that the global phase portrait in the local chart $U_1$ on the infinite sphere is shown
in Figure \ref{u1}: system (\ref{nu11}) has only one equilibrium $(z_1,z_2)=(0,0)$, which is a center,
two invariant straight lines
$$
L_1:\{(z_1,z_2)|z_2=\frac{1}{\sqrt{A}}\},~~L_2:\{(z_1,z_2)|z_2=-\frac{1}{\sqrt{A}}\}
$$
and the region between the two lines are fulfilled with periodic cycles.

\begin{figure}
  \centering
  \includegraphics[width=230pt]{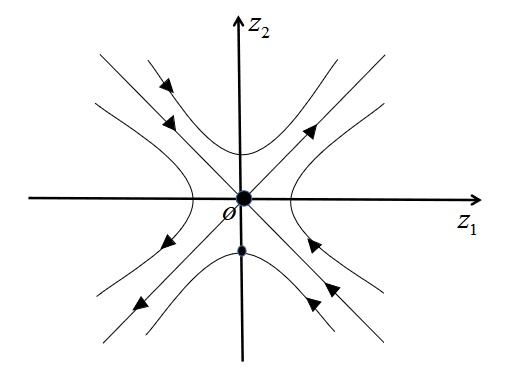}
  \caption{Dynamics of the GD model on the Poincar\'{e} sphere at infinity in the local chart $U_2$}
  \label{u2}
\end{figure}
\emph{In the local charts $U_2$}  We make the change of variables $(x,y,z)=(z_1z_3^{-1},z_3^{-1},z_2z_3^{-1})$
   and  ${\rm d}{\tau}=z_3^{-1}{\rm d}t$, which yields the Poincar\'{e} compactification $p(X)$ of the GD model in the local chart $U_{2}$
\begin{equation}\label{nu2}
\begin{cases}
\frac{\textrm{d}z_{1}}{\textrm{d}{\tau}}=-R_az_1z_3^2+Cz_2z_3-\sigma z_1z_3+z_1^2z_2+Az_2+z_1z_3\\
\frac{\textrm{d}z_{2}}{\textrm{d}{\tau}}=-R_az_2z_3^2+z_1z_2^2+z_1\\
\frac{\textrm{d}z_{3}}{\textrm{d}{\tau}}=-z_3(R_az_3^2-z_1z_2-z_3).
\end{cases}
\end{equation}
For $z_{3}=0$, system (\ref{nu2}) reduces to
\begin{equation}\label{nu22}
\frac{\textrm{d}z_{1}}{\textrm{d}{\tau}}=z_1^2z_2+Az_2,~~~~~
\frac{\textrm{d}z_{2}}{\textrm{d}{\tau}}=z_1z_2^2+z_1.
\end{equation}
By some direct calculations, we find that (\ref{nu22}) is also integrable with a rational first integral
$$
\Phi(z_1,z_2)=\frac{z_2^2+1}{z_1^2+A}.
$$
Then, one can easily get the global phase portrait of system (\ref{nu22}) is shown
in Figure \ref{u2}: system (\ref{nu22}) has only one equilibrium $(z_1,z_2)=(0,0)$,
which is a saddle.

\begin{figure}
  \centering
  \includegraphics[width=230pt]{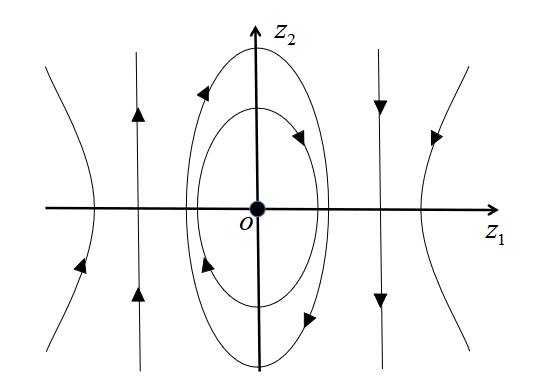}
  \caption{Dynamics of the GD model on the Poincar\'{e} sphere at infinity in the local chart $U_3$}
  \label{u3}
\end{figure}
\emph{In the local charts $U_3$}
Taking the transformation  $(x,y,z)=(z_1z_3^{-1},z_2z_3^{-1},z_3^{-1})$ and ${\rm d}\tau=z_3^{-1}{\rm d}t$,
we obtain the Poincar\'{e} compactification $p(X)$ of the GD model in the local chart $U_{3}$
 \begin{equation}\label{nu3}
\begin{cases}
\frac{\textrm{d}z_{1}}{\textrm{d}{\tau}}=-\sigma z_1z_3-z_1^2z_2+Az_2+Cz_3+z_1z_3\\
\frac{\textrm{d}z_{2}}{\textrm{d}{\tau}}=R_az_3^2-z_1z_2^2-z_1\\
\frac{\textrm{d}z_{3}}{\textrm{d}{\tau}}=-z_3(z_1z_2-z_3).
\end{cases}
\end{equation}
Setting $z_{3}=0$,
system (\ref{nu3}) restricted to the $z_{1}z_{2}-$plane reads
\begin{equation}\label{nu33}
\frac{\textrm{d}z_{1}}{\textrm{d}{\tau}}=-z_1^2z_2+Az_2,~~
\frac{\textrm{d}z_{2}}{\textrm{d}{\tau}}=-z_1z_2^2-z_1.
\end{equation}
Clearly, system (\ref{nu33}) has a rational first integral
$$
\Phi(z_1,z_2)=\frac{z_2^2+1}{-z_1^2+A}.
$$
Proceeding as above, we see that  the global phase portrait of (\ref{nu33}) in the local
chart $U_3$ on the infinite sphere is shown in Figure \ref{u3}:
system (\ref{nu11}) has only one equilibrium $(z_1,z_2)=(0,0)$,
which is a center,
two invariant straight lines
$$
\widehat L_1:\{(z_1,z_2)|z_1=\frac{1}{\sqrt{A}}\},~~\widehat L_2:\{(z_1,z_2)|z_1=-\frac{1}{\sqrt{A}}\}
$$
and the region between the two lines are
fulfilled with periodic cycles.

\begin{figure}
  \centering
  \includegraphics[width=230pt]{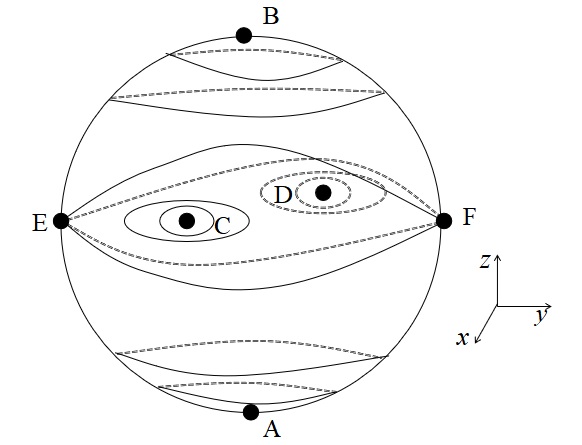}
  \caption{Dynamics of the GD model on the sphere at infinity.}
  \label{u123}
\end{figure}

Finally,
we mention that
the flow in the local charts $V_{i},~i=1,2,3$ is the same as the flow in the local charts $U_{i},~i=1,2,3$ reversing the time,
because the compactified vector field $p(X)$ in $V_{i}$ coincides with the vector field $p(X)$ in $U_{i}$ multiplied by $(-1)^{n-1}$,
where $n$ is the dimension of the considered system.
Then, combining with the above analysis and Figure \ref{u123},
we can obtain the global structure of the dynamical behavior of the GD model on the Poincar\'{e} sphere at infinity as
follows:

\begin{thm}
The phase portrait of
the GD model on the Poincar\'{e} sphere is shown in Figure \ref{u123}: it has four centers $A,B,C,D$
 localized at the endpoints
of the $x$-axis and $y$-axis, and has two saddles $E,F$ localized at the endpoints
of the $z$-axis.
\end{thm}

\section{Conclusion}

In this study, the Darboux integrability,
local bifurcations and global dynamics of a three-mode model describing the convection of viscous incompressible fluid
motion inside the ellipsoid are examined.

In section 3,
we show that the GD model is not Darboux integrable for any values of the physical parameters $A,\sigma,C, R_a$.
This coincides with the numerical fact \cite{gb0,gb1,gb2} that the GD model has chaotic behaviors for a large range of its parameters.
Theorem \ref{th} also implies that any search for a closed-form, analytical solution
for the GD model is bound to fail.

In section 4,
we show that the GD model may undergo two different transitions when the control parameter $R_a$ increases for fixed $A,C,\sigma$.
For the first transition at $R_a=R_a^*$,
the GD model bifurcates from an unstable steady state $E_0$ to two stable steady states $E_{\pm}$.
Interestingly, as $R_a$ further increases,
the GD model exhibits a second transition at $R_a=R_a^{**}$ and bifurcates from the steady states $E_{\pm}$ to periodic solutions in the case of $L(A,\sigma,C)>0$,
whereas no new transition occurs in the case of  $L(A,\sigma,C)\leq0$.

In section 5,
we provide a complete description of the global dynamics of the GD model at infinity.
Our result shows that the parameters $C,R_a,\sigma$
do not effect its global dynamics at infinity, and the
dynamics at the infinity for different positive values of the parameter $A$ are topologically equivalent
although it depends on the parameter $A$.
This means that the parameter $A$ just yields quantitative,
but not qualitative changes for the dynamics at infinity of the GD model.

The theoretical
analysis and  numerical observations in this paper may be useful
both in mathematical and physical research areas. In the following research, more in-depth discussions and research results will be provided.

\appendix

\section{Roots of cubic polynomials}
For the convenience, we recall some basic results used in this work.
Consider a general cubic polynomial
\begin{align}\label{cub}
\lambda^3+p_1\lambda^2+p_2\lambda+p_3=0,
\end{align}
where $p_i,~i=1,2,3$ are real numbers.

The next result is the well-known Routh-Hurwitz criterion in $\mathbb{R}^3$.
\begin{prop}\cite{rh}
The real parts of all the roots $\lambda$ of (\ref{cub}) are negative if and only if
$$
p_1>0,~p_2>0,~p_3>0,~p_1p_2>p_3.
$$
\end{prop}

We say that a root $\lambda$ of (\ref{cub}) is critical if the real part of $\lambda$ is zero. The following result is elementary but useful in the analysis of bifurcations of differential systems in $\mathbb{R}^3$.
\begin{prop}\label{ap1}
The following statements hold for (\ref{cub}).

(1) It has a simple zero root and no other critical roots if and only if $p_3=0,~p_2<0$.

(2) It has a simple pair of purely imaginary roots $\pm iw_0 (w_0>0)$ and no other critical roots if and only if $p_1p_2=p_3,~p_2=w_0^2>0,~p_3\neq 0$.

(3) It has a zero root of (algebraic) multiplicity two and no other critical roots if and only if $p_2=p_3=0,~p_1\neq0$.

(4) It has a simple zero root and a simple pair of purely imaginary roots $\pm iw_0 (w_0>0)$ if and only if $p_2=p_3=0$
and $p_1=w_0^2>0$

\end{prop}
Since its proof is straightforward, we omit it.

\section{Darboux theory of integrability}

To prove Theorem \ref{th}, we need the following results.
\begin{prop}\label{PP1}
Suppose that a polynomial vector field $X$ defined in $\mathbb{R}^n$ of degree
$m$ admits $p$ Darboux polynomials $f_{i}$ with cofactor $K_{i}$ for $i=1,\cdots,p,$
and $q$ exponential factors $E_{j}=\exp(g_{j}/h_{j})$ with cofactors $L_{j}$ for
$j=1,\cdots,q$. If there exist $\lambda_{i}$, $\mu_{j}\in\mathbb{R}$ not all zero
such that
\begin{equation*}\label{P5}
\sum_{i=1}^{p}\lambda_{i}K_{i}+\sum_{j=1}^{q}\mu_{j}L_{j}=0,
\end{equation*}
then the following real (multivalued) function of Darboux type
\begin{equation*}\label{P6}
f_{1}^{\lambda_{1}}\cdots f_{p}^{\lambda_{p}}E_{1}^{\mu_{1}}\cdots E_{p}^{\mu_{q}},
\end{equation*}
substituting $f_{i}^{\lambda_{i}}$ by $|f_{i}|^{\lambda_{i}}$ if $\lambda_{i}\in\mathbb{R}$,
is a first integral of the vector field $X$.
\end{prop}

The proof can be seen in \cite{dar1}.

\begin{prop}\label{PP2}
The following statements hold.\\
(a) If $e^{g/h}$ is an exponential factor for the polynomial differential system (\ref{I1})
and $h$ is not a constant polynomial, then $h$ is a Darboux polynomial\\
(b) Eventually $e^{g}$ can be an exponential factor, coming from the multiplicity of the
infinite invariant plane.
\end{prop}

For a proof of this result see \cite{dar2}.

\begin{prop}\label{PP3}
Let $f$ be a polynomial and $f=\prod_{j=1}^{s}f_{j}^{\alpha_{j}}$ be its decomposition into irreducible
factors in $\mathbb{R}[x,y,z]$. Then $f$ is a Darboux polynomial of system (\ref{I1}) if and only if all
the $f_{j}$ are Darboux polynomials of system (\ref{I1}). Moreover, if $K$ and $K_{j}$ are the cofactors of
$f$ and $f_{j}$, then $K=\sum_{j=1}^{s}\alpha_{j}K_{j}$.
\end{prop}

The proof of Proposition \ref{PP3} can be found in \cite{dar3}.

\end{document}